\documentclass[11pt]{article}

\usepackage{a4wide}
\usepackage{theorem}
\usepackage{amsmath}
\usepackage{array}
\usepackage{amssymb}
\usepackage{amsfonts}
\usepackage[english]{babel}
\usepackage{epsf}
\usepackage{epsfig}
\usepackage{graphicx}

\newcommand{\R}{\mathbb{R}}
\newcommand{\N}{\mathbb{N}}

\newcommand{\C}{\mathbb{C}}
\newcommand{\E}{\mathbb{E}}
\newcommand{\Z}{\mathbb{Z}}

\newcommand{\Sd}{\textup{Sd}} 
\newcommand{\Image}{\textup{Im}} 
 
\newcommand{\Lk}{\textup{Lk}} 
 
\newcommand{\dvol}{d\textup{vol}}

\theoremstyle{definition}
\newtheorem{thm}{Theorem}
\newtheorem{lem}[thm]{Lemma}
\newtheorem{prop}[thm]{Proposition}
\newtheorem{cor}[thm]{Corollary}
\newtheorem{rem}[thm]{Remark}
\newtheorem{ex}[thm]{Example}
\newtheorem{defn}[thm]{Definition}
\newenvironment{pf}{{\bf{Proof.}}}{\hfill $\Box$\\}

\begin{document}

\title{Asymptotic measures and links in simplicial complexes}

\author{Nerm{\accent95\i}n Salepc{\accent95\i} and Jean-Yves Welschinger}

\maketitle

\begin{abstract}
We introduce canonical measures on a locally finite simplicial complex $K$ and study their asymptotic behavior under infinitely many barycentric subdivisions. We also compute the face polynomial of the asymptotic link and dual block of a simplex in the $d^{th}$ barycentric subdivision $\Sd^d(K)$ of $K$, $d\gg0$. It is almost  everywhere  constant. When $K$ is finite, we study the limit face polynomial of $\Sd^d(K)$ after F.~Brenti-V.~Welker and E.~Delucchi-A.~Pixton-L.~Sabalka.

\vspace{0.5cm}


\end{abstract}

\section{Introduction}
Let  $K$ be a finite $n$-dimensional simplicial complex and $\Sd^d(K), d\geq0$, be its successive barycentric subdivisions, see~\cite{M}. We denote by $f_p(K), p\in\{0,\ldots,n\}$, the \emph{face number} of $K$, that is the number of $p$-dimensional simplices of $K$ and by $q_K(T)=\sum_{p=0}^{n}f_p(K)T^p$ its \emph{face polynomial}. The asymptotic of $f_p^d(K)=f_p(\Sd^d(K))$ has been studied in~\cite{BW} and~\cite{DPS}, it is equivalent to $q_{p,n}f_n(K)(n+1)!^d$ as $d$ grows to $+\infty$, where $q_{p,n}>0$. Moreover, it has been proved in~\cite{BW} that the roots of the \emph{limit face polynomial} $q_n^{\infty}(T)=\sum_{p=0}^{n}q_{p,n}T^p$ are all simple and real in $[-1,0]$, and  in~{\cite{DPS} that this polynomial is symmetric with respect to the involution $T\to -T-1$, see Theorem~\ref{Thm_DPS}. We first observe that this symmetry actually follows from a general symmetry phenomenon obtained by I.~G.~Macdonald in~\cite{Mac} which can be formulated as follows, see Theorem~\ref{Thm_$R_K(T)$}.  We set $R_K(T)=Tq_K(T)-\chi(K)T$.

\begin{thm} [Theorem 2.1,~\cite{Mac}] \label{Intro_$R_K(T)$}Let $K$ be a triangulated compact homology $n$-manifold. Then,  $R_K(-1-T)=(-1)^{n+1}R_K(T).$ 
\end{thm}

Recall that a homology $n$-manifold is a topological space $X$ such that  for every $x\in X$,
the relative homology $H_*(X,X\setminus\{x\};\Z)$ is isomorphic to  $H_*(\R^n,\R^n\setminus\{0\}; \Z)$. Any smooth or topological manifold is thus a homology manifold and Poincar\'e duality holds true in such compact homology manifolds, see~\cite{M}. 

We also observe the following theorem (see Corollary~\ref{Cor_$R_K(T)$}  and Theorem~\ref{Thm_onlyroot}), the first part of which is a corollary of Theorem~\ref{Intro_$R_K(T)$}  which has been independently (not as a corollary of Theorem~\ref{Intro_$R_K(T)$}) observed by T. Akita~\cite{A}.

\begin{thm}\label{Intro_chiK}
Let $K$ be a compact triangulated homology manifold of even dimension. Then $\chi(K)=q_K(-\frac{1}{2}).$ 

Moreover, $t=-1$ in odd dimensions and $t=-1$ together with  $t=-\frac{1}{2}$ in even dimensions are the only complex values $t$ on which $q_K(T)$ equals $\chi(K)$ for every compact triangulated homology manifold of the given dimension.
\end{thm}

Having spheres in mind for instance, Theorem~\ref{Intro_$R_K(T)$} and  Theorem~\ref{Intro_chiK} exhibit a striking behavior of simplicial structures compared to cellular structures. In~\cite{SW2}, we  also provide a probabilistic proof of the first part of Theorem~\ref{Intro_chiK}.

The limit face polynomial $q_n^\infty(T)$ remains puzzling, but we have been able to prove the following  result, see   Proposition~\ref{Prop_qpn} and Corollaries~\ref{Cor_$L_i$} and~\ref{Cor_Eigenvector}. 

Let $L_j(T)=\frac{1}{j!}\prod\limits_{i=0}^{j-1}(T-i) \in \R[T], j\geq1$, be the $j^{th}$ Lagrange polynomial and set $L_0=1$.

\begin{thm}
Let $\Lambda^t=(\lambda_{j,i})$ be the upper triangular matrix of the vector $(T^j)_{j\geq0}$ in the base $(L_i)_{i\geq 0}$. Then, $(q_{p,n})_{0\leq p\leq n}$
is the eigenvectors of $\Lambda^t$ associated to the eigenvalue $(n+1)!$ normalized in such a way that $q_{n,n}=1$ and $q_{p,n}=0$ if $p>n$. Moreover  for every $0\leq p<n,$
$$q_{p,n}=\sum\limits_{(p_1,\ldots, p_j)\in\mathcal{P}_{p,n}}\frac{ \lambda_{n+1,p_{j}}\lambda_{p_{j},p_{j-1}}\dots \lambda_{p_2,p_1}}{(\lambda_{n+1,n+1}-\lambda_{p_{j}, p_{j}})\ldots(\lambda_{n+1,n+1}-\lambda_{p_1,p_1})},$$
where $\mathcal{P}_{p,n}=\{(p_1,\ldots,p_j)\in \N^j | j\geq1 \mbox{ and } p+1=p_1<\ldots<p_j<n+1\}$.
\end{thm}

Our main purpose in this paper is to refine this asymptotic study of the face polynomial by introducing a canonical measure on $\Sd^d(K)$ and study the density of links in $\Sd^d(K)$ with respect to these measures. For every $0\leq p\leq n$, set $\gamma_{p,K}=
\sum_{\sigma\in K^{[p]}}\delta_{\hat{\sigma}}$, where $\delta_{\hat{\sigma}}$
denotes the Dirac measure on the barycenter $\hat{\sigma}$ of $\sigma$ and $K^{[p]}$ the set of $p$-dimensional simplices of $K$. Likewise, for every $d\geq 0, $ we set $\gamma_{p,K}^d=\frac{1}{(n+1)!^d}\gamma_{p, \Sd^d(K)}$, which provides a canonical sequence of Radon measures on the underlying topological space $|K|.$ The latter is also equipped with the measure $\dvol_K=\sum_{\sigma\in K^{[n]}}(f_{\sigma})_* \dvol_{\Delta_n},$ where $f_\sigma :\Delta_n\to \sigma$ denotes a simplicial isomorphism between the standard $n$-simplex $\Delta_n$ and the simplex $\sigma$, and $\dvol_{\Delta_n}$  denotes the Lebesgue measure normalized in such a way that $\Delta_n$ has volume 1, see Section~\ref{Sect_CanMes}.
We prove the  following, see Theorem~\ref{Thm_Gammapd}.

\begin{thm}\label{Intro_Gammapd}
For every $n$-dimensional  locally finite simplicial complex $K$ and every $0\leq p\leq n$, the measure $\gamma_{p,K}^d$
weakly converges to $q_{p,n}\dvol_K$ as $d$ grows to $+\infty$.
\end{thm}

When $K$ is finite, Theorem~\ref{Intro_Gammapd} recovers the asymptotic of $f_p^d(K)$ as $d$ grows to $+\infty$, by integration of the constant function 1. 
Recall that the \emph{link of a simplex} $\sigma$ in $K$ is by definition
$\Lk(\sigma,K)=\{\tau\in K| \sigma \mbox{ and } \tau \mbox{ are disjoint and both are faces of a simplex in } $K$\}.$ Likewise, the \emph{block dual} to $\sigma$ is the set $D(\sigma)=\{[\hat{\sigma}_0,\ldots, \hat{\sigma}_p] \in \Sd(K)| p\in\{0,\ldots,n\} \mbox{ and } \sigma_0=\sigma\}$, see~\cite{M}. Recall that the simplices of $\Sd(K)$ are by definition of the form $[\hat{\sigma}_0,\ldots, \hat{\sigma}_p]$, where $\sigma_0<\ldots<\sigma_p$ are simplices of $K$ with $<$ meaning  being a proper face. The dual blocks form a partition of $\Sd(K),$ see~\cite{M}, and the links $\Lk(\sigma,K)$ encodes in a sense the local complexity of $K$ near $\sigma.$ We finally prove the following, see Theorem~\ref{Thm_qLk} and Theorem~\ref{Thm_qDsigma}.

\begin{thm}\label{Intro_qLk}
For every $n$-dimensional locally finite  simplicial complex $K$ and every $0\leq p<n,$ the measure 
$q_{\Lk(\sigma,\Sd^d(K))}(T)d\gamma_{p,K}^d(\sigma)$  (with value in $\R_{n-p-1}[T]$) weakly converges to $\left(\sum_{l=0}^{n-p-1}q_{p+l+1,n} f_p(\Delta_{p+l+1})T^l\right)\dvol_K$ as $d$ grows to $+\infty.$
\end{thm}

And likewise 
\begin{thm}\label{Intro_qDsigma}
For every $n$-dimensional locally finite  simplicial complex $K$ and every $0\leq p\leq n$, the measure $q_{D(\sigma)}(T)d\gamma^d_{p,K}(\sigma)$ weakly converges to  $\sum_{l=0}^{n-p}\left(\sum_{h=l}^{n-p}q_{p+h,n}f_p{(\Delta_{p+h})}\lambda_{h,l}\right)T^l \dvol_K$ as $d$ grows to $+\infty.$
\end{thm}

From these theorems we see that asymptotically, the complexity of the link and the dual block is almost everywhere constant with respect to $\dvol_K$. In~\cite{SW2}, we study the asymptotic topology of a random subcomplex in a finite simplicial complex $K$ and its successive barycentric subdivisions. It turns out that the Betti numbers of such a subcomplex get controlled by the measures given in Theorem~\ref{Intro_qDsigma}.

\vspace{0.5 cm}

\textbf{Acknowledgement :}
The second author is partially supported by the ANR project MICROLOCAL (ANR-15CE40-0007-01).

\section{The face polynomial of a simplicial complex} \label{Sect_FaceP}

\subsection{The symmetry property}\label{Ssect_Sym}
Let $K$ be a finite $n$-dimensional simplicial complex. We set $ R_K(T)=Tq_K(T)-\chi(K)T$, where $q_K(T)=\sum_{p=0}^{n}f_p(K)T^p$ and $\chi (K)$ is the Euler characteristic of $K$, so that $R_K(0)=R_K(-1)=0.$

\begin{ex} \begin{enumerate}
\item If $K=\partial\Delta_{n+1},$ then $Tq_K(T)=(1+T)^{n+2}-1-T^{n+2}.$
\item If $K=S^0\ast\ldots\ast S^0$ is the $n^{th}$ iterated suspension of the $0$-dimensional sphere, then $
R_K(T)=Tq_K(T)-T\chi(K)=\begin{cases} 
(2T+1)((2T+1)^n-1) & \mbox{if $n$ is even},\\
 (2T+1)^{n+1}-1& \mbox{if $n$ is odd}.\\
  \end{cases}$
\end{enumerate}
\end{ex}

Recall that if $K$ is a triangulated compact homology $n$-manifold, 
its face numbers satisfy the following  Dehn-Sommerville relations (\cite{Klee},  see also for example~\cite{Klain}):
$$\forall\, 0\leq p\leq n,\, \, f_p(K)=\sum\limits_{i=p}^{n}(-1)^{i+n}\binom{i+1}{p+1}f_i(K).$$

The Dehn-Sommerville relations imply that $R_K(T)$ satisfy the following  striking symmetry property observed by  I.G.~Macdonald~\cite{Mac} which we recall here together with a proof for the reader's convenience.

\begin{thm} [Theorem 2.1,~\cite{Mac}] \label{Thm_$R_K(T)$}Let $K$ be a triangulated compact homology $n$-manifold. Then,  $R_K(-1-T)=(-1)^{n+1}R_K(T).$ 
\end{thm}

\begin{pf} Observe that 
$$\begin{array}{rcl}
R_K(-1-T)&=&\sum\limits_{p=0}^{n}f_p(K)(-1-T)^{p+1}+\chi(K)(1+T)\\
&=&\sum\limits_{p=0}^{n}f_p(K)(-1)^{p+1}\sum\limits_{q=0}^{p+1}\binom{p+1}{q} T^q+\chi(K)(1+T)\\
&=&\sum\limits_{p=0}^{n}f_p(K)(-1)^{p+1}\sum\limits_{q=0}^{p}\binom{p+1}{q+1}T^{q+1}+\chi(K)T\\
&=&\sum\limits_{q=0}^{n}T^{q+1}\sum\limits_{p=q}^n\binom{p+1}{q+1}f_p(K)(-1)^{p+1}+\chi(K)T.\\
\end{array}$$

Then, the Dehn-Sommerville relations imply
$$\begin{array}{rcl}
R_K(-1-T)&=&-\sum\limits_{q=0}^{n}T^{q+1}(-1)^nf_q(K)+\chi(K)T\\
&=&(-1)^{n+1}R_K(T)+(1+(-1)^{n+1})\chi(K)T.\end{array}$$
Now, if $n$ is even, $1+(-1)^{n+1}=0$ while if $n$ is odd, $\chi(K)=0$ by Poincar\'e duality with $\Z/2\Z$ coefficients, see~\cite{M}. In both cases, we get $R_K(-1-T)=(-1)^{n+1}R_K(T).$
\end{pf}

\begin{cor} \label{Cor_$R_K(T)$} Let $K$ be a triangulated compact homology $n$-manifold. 
\begin{enumerate} 
\item If $n$ is even, then $q_K(-\frac{1}{2})=\chi(K)$.
\item If $n$ is odd, the polynomial $Tq_K(T)$ is preserved by the involution $T\to -1-T$.
\item If $\chi(K)\leq 0$,  the real roots of $R_K(T)=Tq_K(T)-\chi(K)T$ lie on the interval $[-1, 0].$
\end{enumerate}
\end{cor}

\begin{pf}
When $n$ is even,  $R_K$ has an odd number of real roots, invariant under the involution $T\to -1-T$ whose unique fixed point is $-\frac{1}{2}$. Theorem~\ref{Thm_$R_K(T)$} thus implies that $R_K(-\frac{1}{2})=0$. Hence the first part. 
When $n$ is odd, $\chi(K)=0$ by Poincar\'e duality so that $R_K(T)=Tq_K(T)$ and the second part.
Finally, if $\chi(K)\leq 0$, the coefficients of the polynomial $R_K(T)$ are all positive, so that its real roots are all negative. It thus follows from  Theorem~\ref{Thm_$R_K(T)$} that they lie on the interval $[-1,0]$.
\end{pf}

\begin{rem} 
  The first part of Corollary~\ref{Cor_$R_K(T)$} was independently (not as a corollary of Theorem~\ref{Thm_$R_K(T)$}) observed by T.~Akita~\cite{A}. In~\cite{SW2}, we provide a probabilistic proof of it. 
  
The third part of Corollary~\ref{Cor_$R_K(T)$} always holds true when $n$ is odd, since then $\chi(K)=0$.
\end{rem}

The first part of Corollary~\ref{Cor_$R_K(T)$} raise the following question: given some dimension $n$, what are the 
universal parameters $t$ such that $q_K(t)=\chi(K)$ for every compact triangulated homology $n$-manifolds?
We checked that $t=-1$ in odd dimensions and $t=-1$ with $t=-\frac{1}{2}$ in even dimensions   are the  only ones, see Theorem~\ref{Thm_onlyroot}.

\subsection{The asymptotic face polynomial}\label{Ssect_Asymp}
Let  $f(K)=(f_0(K),f_1(K),\ldots,f_n(K))$ be the \emph{face vector} of $K$, that is the vector formed by the face numbers of the finite simplicial complex $K$.
Now,  for every $d>0$, we set $f_p^d(K)=f_p(\Sd^d(K))$, where $\Sd^d(K)$ denotes the $d^{th}$ barycentric subdivision of $K$.  How does the face vector change under barycentric subdivisions and what is the asymptotic behavior of $f^d(K)=(f_0^d(K),f_1^d(K) \ldots, f_n^d(K))$? These questions have been treated in~\cite{BW},~\cite{DPS}, leading to the following. 

\begin{thm} [\cite{BW},~\cite{DPS}] \label{Thm_DPS}   For every $0\leq p\leq n$, there exist $q_{p,n}>0$ such that for every  $n$-dimensional finite simplicial complex $K,$ $\lim_{d\to +\infty} \frac{f_p^d(K)}{(n+1)!^df_n(K)}=q_{p, n}.$   

Moreover, the $n+1$ roots of the polynomial $Tq_n^\infty(T)$ are simple, belong to the interval $[-1,0]$ and are symmetric with respect to the involution $T\in \R\mapsto -T-1\in \R$ whenever $n>0$, where $q_n^\infty(T)=\sum\limits_{p=0}^{n}q_{p, n} T^p$.
\hfill$\square$\\
\end{thm} 

The symmetry property of $Tq_n^\infty(T)$ follows from Theorem~\ref{Thm_$R_K(T)$} and the first part of Theorem~\ref{Thm_DPS}, since the Euler characteristic remains unchanged under subdivisions. This symmetry has been observed in~\cite{DPS} (with a different proof).
It implies that  $q_n^\infty(-1)=0$ and that  $q_n^\infty(-\frac{1}{2})=0$ whenever $n$ is even, as the number of roots of  $Tq_n^\infty(T)$ is then odd and $-\frac{1}{2}$ is the unique fixed point of the involution. 

\begin{thm}\label{Thm_onlyroot}
The reals $t=-1$ if $n$ is odd and $t= -1$ together with $t=-\frac{1}{2}$ if $n$ is even are the only complex values on which the face polynomial $q_K(T)=\sum_{p=0}^{\dim K}f_p(K)T^p$ equals $\chi(K)$ for every  compact triangulated homology $n$-manifold $K$.
\end{thm}

\begin{pf}
Let us equip the $n$-dimensional sphere with the triangulation given by the boundary of the $(n+1)$-simplex $\Delta_{n+1}$. Then, for every $0\leq p\leq  n$, $f_p(S^n)=\binom{n+2}{p+1}$  and $q_{S^n}(T)=\frac{1}{T}\big((1+T)^{n+2}-1-T^{n+2}\big).$ Now, the polynomial $q_{S^n}(T)-\chi(S^n)$ has only one real root if $n$ is odd
 and two real roots if $n$ is even. Indeed, differentiating the polynomial $Tq_{S^n}(T)-\chi(S^n)T$ once if $n$ is odd and twice if $n$ is even, we get, up to a factor, $(1+T)^{n+1}-T^{n+1}$  or respectively $(1+T)^{n}-T^{n}$ which vanishes only for $t=-\frac{1}{2}$ on the real line. From Rolle's theorem we deduce that 0 and -1 (respectively $0, -\frac{1}{2}, -1$) are the only real roots of $Tq_{S^n}(T)-\chi(S^n)T$ when $n$ is odd (respectively, when $n$ is even). 
 
 Finally,  if $t_0\in \C$ is such that $q_K(t_0)=\chi(K)$ for all triangulated manifolds of a given dimension $n$, then in particular, $R_{\Sd^d(K)}(t_0)=0$ for every $d>0$.
 Dividing  by $f_n(K)(n+1)!^d$ and passing to the limit, we deduce that $q_n^\infty(t_0)=0$. But from Theorem~\ref{Thm_DPS}
 we know that the roots of $Tq^\infty(T)$ are all real, hence the result.
\end{pf}

Let now $\Lambda=(\lambda_{i,j})_{i, j\geq1}$ be the infinite lower triangular matrix whose entries $\lambda_{i,j}$ are
the numbers of  interior $(j-1)$-faces  on the subdivided standard simplex $\Sd(\Delta_{i-1})$ and let $\Lambda_n=(\lambda_{i,j})_{1\leq i,j\leq n+1}$, see Figure~\ref{Lambda}.
The diagonal entries of $\Lambda$ are given by Lemma~\ref{Lem_lambdaij}. We set as a convention $\lambda_{0,0}=1$ and $\lambda_{l,0}=0$ whenever $l>0$. 

\begin{figure}[h]
   \begin{center}
    \includegraphics[scale=0.3]{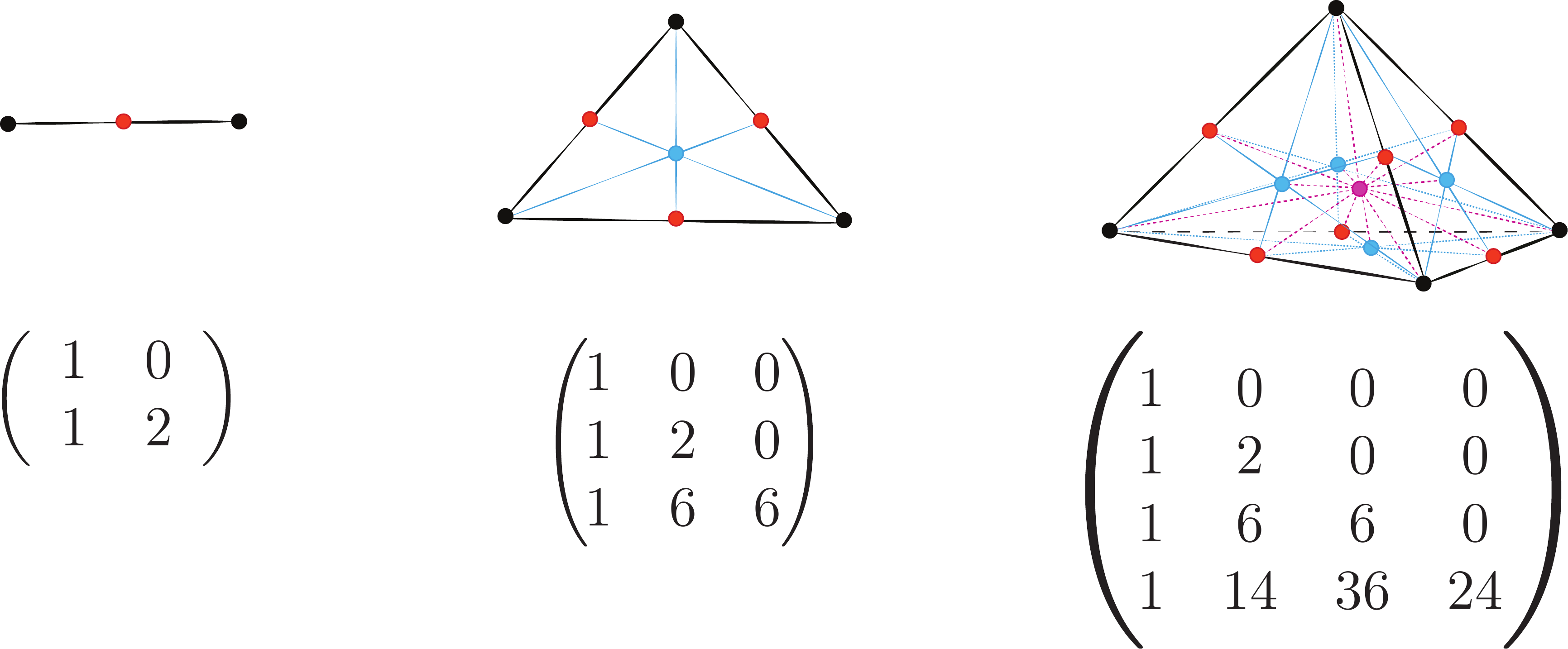}
    \caption{The matrix $\Lambda_n$ for $n=1, 2, 3$.}
    \label{Lambda}
      \end{center}
 \end{figure}

\begin{lem}\label{Lem_lambdaij} For every $1\leq j\leq i$,  $\lambda_{i,j}=\sum_{p=j-1}^{i-1}\binom{i}{p}\lambda_{p,j-1}$ where $\binom{i}{j}$ denotes the binomial coefficient. In particular, $\lambda_{i,i}=i!$. 
\end{lem}
\begin{pf} 
The interior $(j-1)$-faces of $\Sd(\Delta_{i-1})$ are cones over the $(j-2)$-faces of the boundary of $\Sd(\Delta_{i-1})$. The latter are interior to some $(p-1)$-simplex of $\partial\Delta_{i-1}$, $j-1\leq p\leq i-1$. The result follows from the fact that for every $1\leq p\leq i-1$, $\partial \Delta_{i-1}$ has $\binom{i}{p}$ many $(p-1)$-dimensional faces while each such face contains $\lambda_{p,j-1}$ many $(j-2)$-dimensional faces of $\Sd(\Delta_{i-1})$ in its interior.
\end{pf}

The first part of Theorem~\ref{Thm_DPS} is basically deduced  in~\cite{BW},~\cite{DPS} from the following observation: for every $n$-dimensional  finite simplicial complex $K$, the face vector $f(\Sd(K))$
is deduced from the face vector $f(K)$  by multiplication on the right by $\Lambda_n$, that is $f(\Sd(K))=f(K)\Lambda_n$, while the matrix $\Lambda_n$ is diagonalizable with eigenvalues given by Lemma~\ref{Lem_lambdaij}. 
  
We deduce from~\cite{BW},~\cite{DPS} that the vector $(q_{p,n})_{0\leq p\leq n}$ is the eigenvector of $\Lambda_n^t$ associated to the eigenvalue $\lambda_{n+1,n+1}=(n+1)!$ normalized by the relation $q_{n,n}=1$. A geometric proof of this fact will be given in Section~\ref{Sect_Links}, see Corollary~\ref{Cor_Eigenvector}. This observation makes it possible to compute $q_{p,n}$ in terms of the coefficients $\lambda_{i,j}$.

\begin{prop}\label{Prop_qpn} Let $0\leq p< n$ and let $\mathcal{P}_{p,n}=\{(p_1,\ldots,p_j)\in \N^j | j\geq1 \mbox{ and } p+1=p_1<\ldots<p_j<n+1\}$. 
Then

$$q_{p,n}=\sum\limits_{(p_1,\ldots, p_j)\in\mathcal{P}_{p,n}}\frac{ \lambda_{n+1,p_{j}}\lambda_{p_{j},p_{j-1}}\dots \lambda_{p_2,p_1}}{(\lambda_{n+1,n+1}-\lambda_{p_{j}, p_{j}})\ldots(\lambda_{n+1,n+1}-\lambda_{p_1,p_1})}.$$
\end{prop}
 
 \begin{pf}
Having in mind that $\Lambda_n$ is a  lower triangular matrix and by Lemma~\ref{Lem_lambdaij}, $(n+1)!=\lambda_{n+1,n+1}$. The equation $\Lambda_n^t (q_{p,n})=(n+1)!(q_{p,n})$ results in the following system.

For all $0\leq p<n$,
$$q_{p,n}=\sum\limits_{k=0}^{n-p-1}\frac{\lambda_{n+1-k, p+1} q_{n-k,n}}{(\lambda_{n+1,n+1}-\lambda_{p+1,p+1})}.$$ 
The solution of this system is obtained by induction on $r=n-p$ by setting $q_{n,n}=1.$ The result follows from the fact that the partitions $(p_1,\ldots, p_j)$ of integers between $p+1$ and $n+1$ such that $p+1=p_1<\ldots<p_j< n+1$ are obtained (except the one with single term $p_1=p+1$) from those $p+1+s=p'_1<\ldots<p'_j< n+1$ for all $1\leq s \leq r$ by setting $p_1=p+1$ and $p_{i+1}=p'_{i}$ for $i\in\{1,\ldots, j\}$.
 \end{pf}

 Note that the coefficients $\lambda_{i,j}$ of $\Lambda$
can be computed. We recall their values obtained in~\cite{DPS} in the following proposition and suggest an alternative proof. 

\begin{prop}  [Lemma~6.1,~\cite{DPS}]  \label{Prop_DPSLem6}For every $1\leq j\leq i$,  $$\lambda_{i,j}=\sum\limits_{p=0}^{j}\binom{j}{p}(-1)^{j-p}p^i.$$ 
\end{prop}
(The left hand side in Lemma~6.1 of~\cite{DPS} should read $\lambda_{i-1,j-1}$ and our $\lambda_{i,j}$ corresponds to  $\lambda_{i-1, j-1}$ in~\cite{DPS}.)

Let $C=(c_{i,j})_{i,j\geq1}$ be the infinite strictly lower triangular matrix such that  $c_{i,j}=\binom{i}{j}$ for $i>j\geq 1$.
Also, for every $r\geq 1$, set $(I+C)^r=(a^r_{i,j})_{i,j\geq1}$.

\begin{lem}\label{Lem_$I+C$} 
 For every $i\geq j$,  $a^r_{i,j}=\binom{i}{j}r^{i-j}$.
\end{lem}
\begin{pf}
We proceed  by induction on $r$.
The statement holds true for $r=1$.  
In the case $r=2$, for every $i\geq j$,
$$\begin{array}{lcl}
a^2_{i, j}&{=}&\sum\limits_{j\leq p\leq i}\binom{i}{p}\binom{p}{j}\\
&=&\frac{i!}{j!(i-j)!}\sum\limits_{p=j}^{i}\frac{(i-j)!}{(p-j)!(i-p)!}\\
&\stackrel{l=p-j}{=}&\binom{i}{j}\sum\limits_{l=0}^{i-j}\binom{i-j}{l}\\
&=&\binom{i}{j}2^{i-j}.
\end{array}$$
The last line follows from the Newton binomial theorem. 
Now, let us suppose that the formula holds true for $r-1$. Then, likewise, 
$$\begin{array}{lcl}
a^r_{i,j}&=&\sum\limits_{j\leq p\leq i}\binom{i}{p}\binom{p}{j}(r-1)^{p-j}\\
&\stackrel{l=p-j}=&\binom{i}{j}\sum\limits_{l=0}^{i-j}\binom{i-j}{l}(r-1)^l\\
&=&\binom{i}{j}r^{i-j}.\\
\end{array}$$
\end{pf}

\textbf{Proof of Proposition~\ref{Prop_DPSLem6}.} We deduce from Lemma~\ref{Lem_lambdaij} that the $p^{th}$ column of  the matrix $\Lambda$ is obtained from the $(p-1)^{th}$ one by multiplication on the left by $C$, so that it  is equal  to $ C^{p-1}v$ where $v$ denotes the first column of $\Lambda$ with 1 on every entry.
Let $C^r=(c^r_{i,j})_{i\geq1, j\geq1}$, then from the relation $C^r=(I+C-I)^r=\sum\limits_{p=0}^{r}\binom{r}{p}(I+C)^p(-1)^{r-p}$, we deduce thanks to Lemma~\ref{Lem_$I+C$} that  for all $r>0$ and all $i\geq j$, 

$$
\begin{array}{lcl}
c^r_{i,j}&=&\sum\limits_{p=0}^{r}\binom{r}{p}\binom{i}{j}p^{i-j}(-1)^{r-p}\\
&=&\binom{i}{j}\sum\limits_{p=0}^{r}\binom{r}{p}p^{i-j}(-1)^{r-p}
\end{array}
$$
 while $c_{i,j}^r=0$ whenever $i\leq j.$
 From the previous observation we now deduce that for all $i\geq r+1$,

$$
\begin{array}{lcl}
\lambda_{i, r+1}&=&\sum\limits_{j=1}^{i-1}\binom{i}{j}\sum\limits_{p=0}^{r}\binom{r}{p}p^{i-j}(-1)^{r-p}\\
&=&\sum\limits_{p=0}^{r} \binom{r}{p}(-1)^{r-p}p^i\sum\limits_{j=1}^{i-1}\binom{i}{j}p^{-j}\\
&=&\sum\limits_{p=0}^{r}\binom{r}{p}(-1)^{r-p}p^i\big((1+\frac{1}{p})^i-1-p^{-i}\big)\\
&=&\sum\limits_{p=0}^{r}\binom{r}{p}(-1)^{r-p} (p+1)^i -\sum\limits_{p=0}^{r}\binom{r}{p}(-1)^{r-p} p^i. \\
\end{array}$$
Now,  we set $l=p+1$ and get

$$
\begin{array}{lcl}
\lambda_{i, r+1}&=&\sum\limits_{l=1}^{r+1}\binom{r}{l-1}(-1)^{r-l+1} l^i -\sum\limits_{p=0}^{r}\binom{r}{p}(-1)^{r-p} p^i \\
&=&(r+1)^i-\sum\limits_{p=1}^{r}\big(\binom{r}{p-1}+\binom{r}{p}\big)(-1)^{r-p}p^i\\
&=&(r+1)^i-\sum\limits_{p=1}^{r}\binom{r+1}{p}(-1)^{r-p}p^i\\
&=&\sum\limits_{p=1}^{r+1}\binom{r+1}{p}(-1)^{r-p+1}p^i.
\end{array}$$
Hence the result. \hfill$\square$\\

Finally, for every $j\geq 1$, let  $L_j(T)=\frac{1}{j!}\prod\limits_{i=0}^{j-1}(T-i) \in \R[T]$ be the $j^{th}$ Lagrange polynomial, so that $L_j(p)=0$ if $0\leq p< j$ and $L_j(p)=\binom{p}{j}$ if $p\geq j$. We deduce the following interpretation of  the transpose matrix $\Lambda^t$. 

\begin{cor}\label{Cor_$L_i$}
For every $j\geq 1$, $T^j=\sum_{i=1}^j\lambda_{j,i}L_i(T)$.  
\end{cor}
Corollary~\ref{Cor_$L_i$} means that $\Lambda^t$ is the matrix of the vectors $(T^j)_{j\geq 0}$ in the basis $(L_i)_{i\geq 0}$ of $\R[T]$, setting $T^0=L_0=1$.

\begin{pf}
Let $i\geq 1$. Then, for every $l\geq i,$

$$\begin{array}{rcl}
\sum\limits_{p=0}^{l}\binom{l}{p}(-1)^{l-p}L_i(p)&=&\sum\limits_{p=i}^{l}\binom{l}{p}(-1)^{l-p}\binom{p}{i}\\
&=&\binom{l}{i}\sum\limits_{p=i}^{l}\binom{l-i}{l-p}(-1)^{l-p}\\
&=&(-1)^{l-i}\binom{l}{i}\sum\limits_{q=0}^{l-i}\binom{l-i}{q}(-1)^q\\
&=&\delta_{li},
\end{array}$$
where $\delta_{li}=0$ if $l\neq i$ and $\delta_{li}=1$ otherwise. This result also holds true for $l\in\{0,\ldots,i-1\}$.  We deduce that for $0\leq l\leq j$, 
$$\sum_{p=0}^{l}\binom{l}{p}(-1)^{l-p}\big( \sum_{i=0}^{j} \lambda_{j,i}L_i(p)\big)=\lambda_{j,l}.$$ The result now follows from Proposition~\ref{Prop_DPSLem6} and  the fact that  a degree $j$ polynomial   is uniquely determined by its values on the 
$j+1$ integers $\{0,\ldots, j\}$, since the above linear combinations for $l\in \{0,\ldots,j\}$ define an invertible triangular matrix.
\end{pf}

\section{Canonical measures on a simplicial complex}\label{Sect_CanMes}

Let us equip the standard $n$-dimensional simplex $\Delta_n$ with the Lebesgue measure $\dvol_{\Delta_n}$ inherited by some affine embedding of $\Delta_n$ in an Euclidian $n$-dimensional space $\E$ in such a way that the total measure of $\Delta_n$ is 1. This measure does not depend on the embedding $\Delta_n \hookrightarrow \E$ for two such embeddings differ by an affine isomorphism which has constant Jacobian 1. 

\begin{defn} For every $n$-dimensional locally finite simplicial complex $K$, we denote by $\dvol_K$ the measure $\sum_{\sigma\in K^{[n]}}(f_\sigma)_* (\dvol_{\Delta_n})$ of $|K|$ where $K^{[n]}$ denotes the set of $n$-dimensional simplices of $K$ and $f_\sigma:\Delta_n\to \sigma$ a simplicial isomorphism.
\end{defn}

If $K$ is a finite $n$-dimensional simplicial complex, the total measure of $|K|$ is thus $f_n(K)$ and its $(n-1)$-skeleton has vanishing measure.
This canonical measure $\dvol_K$ is Radon with respect to the topology of $|K|$.

Now, for every  $p\in \{0,\ldots, n\}$, we set 
$\gamma_{p,K}=\sum_{\sigma\in K^{[p]}}\delta_{\hat{\sigma}}$, where $\delta_{\hat{\sigma}}$ denotes the Dirac measure on the barycenter $\hat{\sigma}$ of $\sigma$. If $K$ is finite, the total measure $\int_{\sigma\in K^{[p]}}1d\gamma_{p,K}(\sigma)$ equals $f_p(K).$ More generally,  for every $d\geq0$, we set $\gamma_{p,K}^d=\frac{1}{(n+1)!^d}\sum_{\sigma\in \Sd^d(K)^{[p]}}\delta_{\hat{\sigma}}.$  

\begin{thm}\label{Thm_Gammapd}
For every $n$-dimensional  locally finite simplicial complex $K$ and every $p\in \{0,\ldots, n\}$, the measure $\gamma_{p,K}^d$
weakly converges to $q_{p,n}\dvol_K$ as $d$ grows to $+\infty$.
\end{thm}

By weak convergence, we mean that for every continuous function $\varphi$ with compact support in $|K|$, 
$\int_K\varphi d\gamma_{p,K}^d  \underset{d \to + \infty}\longrightarrow q_{p,n}\int_K\varphi \dvol_K.$
In order to prove Theorem~\ref{Thm_Gammapd}, we need first the following lemma.

\begin{lem}\label{Lem_theta}
Let $p\in \{0,\ldots, n\}.$ Then for every $l, d\geq 0$, $$\gamma_{p,\Delta_n}^{l+d}=\frac{1}{(n+1)!^l}\sum \limits_{\sigma\in \Sd^l(\Delta_n)^{[n]}} (f_\sigma)_*(\gamma_{p, \Delta_n}^d)-\theta^l_p(d),$$
where $f_\sigma :\Delta_n\to \sigma$ denotes a simplicial isomorphism and the total measure of $\theta_p^l(d)$ converges to zero as $d$ grows to $+\infty$.
\end{lem}

\begin{pf}
In a subdivided $n$-simplex $\Sd^l(\Delta_n),$ every  $p$-simplex $\tau$ is a face of an $n$-simplex and the number of such $n$-simplices is by definition $f_{n-p-1}(\Lk(\tau, \Sd^l(\Delta_n)))$.
Since $\Sd^{l+d}(\Delta_n)=\Sd^d({\Sd^l(\Delta_n)}),$ we deduce that for every $d\geq 0,$

$$
\begin{array}{lcl}
\gamma^{l+d}_{p, \Delta_n}&=&\frac{1}{(n+1)!^l}\sum_{\sigma\in \Sd^l(\Delta_n)^{[n]}} (f_\sigma)_*(\gamma_{p,\Delta_n}^d)-\\&&\frac{1}{(n+1)!^{l+d}}\sum_{\tau\in\Sd^l(\Delta_n)^{(n-1)}}\big(f_{n-\dim\tau-1}(\Lk(\tau, \Sd^l(\Delta_n)))-1\big) \sum_{\alpha\in \Sd^d(\stackrel{\circ}{\tau})^{[p]}} \delta_{\hat{\alpha}},
\end{array}$$
where $\Sd^l(\Delta_n)^{(n-1)}$ denotes the $(n-1)$-skeleton of $\Sd^l(\Delta_n).$

We thus set 

$\theta_p^l(d)=\frac{1}{(n+1)!^{l+d}}\sum_{\tau\in\Sd^l(\Delta_n)^{(n-1)}}\big(f_{n-\dim \tau-1}(\Lk(\tau, \Sd^l(\Delta_n)))-1\big) \sum_{\alpha\in \Sd^d(\stackrel{\circ}{\tau})^{[p]}} \delta_{\hat{\alpha}}.$

The total mass of this measure $\theta_p^l(d)$ satisfies

$\int_{\Delta_n}1 d \theta_p^l(d)\leq \left( \frac{1}{(n+1)!^l}\sup_\tau \big(f_{n-\dim\tau-1}(\Lk(\tau,\Sd^l(\Delta_n)))-1\big) \times \#\Sd^l(\Delta_n)^{(n-1)}\right)\frac{\sup_{\tau}f_p^d(\stackrel{\circ}{\tau})}{(n+1)!^d}.$

Since $\dim \tau <n$,  we know from Theorem~\ref{Thm_DPS} that $\frac{\sup_\tau f_p^d(\stackrel{\circ}{\tau})}{(n+1)!^d}  \underset{d \to + \infty}\longrightarrow 0.$
Hence the result. 
\end{pf}

\textbf{Proof of Theorem~\ref{Thm_Gammapd}.} 
Let us first assume that $K=\Delta_n$ and let $\varphi\in C^0(\Delta_n)$. We set, for every $l,d\geq0$, $R_{l,d}=\int_{\Delta_n}\varphi d\gamma_{p,\Delta_n}^{l+d}-q_{p,n}\int_{\Delta_n}\varphi \dvol_{\Delta_n}$ and deduce from Lemma~\ref{Lem_theta}   
$$\begin{array}{lcl}
R_{l,d}&=&\frac{1}{(n+1)!^l}\sum_{\sigma\in \Sd^l(\Delta_n)^{[n]}}\big(\int_{\Delta_n} f^*_\sigma \varphi d\gamma_{p,\Delta_n}^d-q_{p,n}\int_{\Delta_n}f^*_\sigma \varphi \dvol_{\Delta_n} \big)- \int_{\Delta_n}\varphi d \theta_p^l(d),\end{array}$$ since by definition $(f_\sigma)_*\dvol_{\Delta_n}=(n+1)!^l\dvol_{\Delta_n}|_\sigma.$ Thus, 

$\begin{array}{lcl}
R_{l,d}&=&\frac{1}{(n+1)!^l}\sum_{\sigma\in \Sd^l(\Delta_n)^{[n]}}\left( \int_{\Delta_n}\big(f^*_\sigma\varphi -\varphi(\hat{\sigma})\big) d\gamma_{p, \Delta_n}^d -q_{p,n}\int_{\Delta_n}\big(f^*_\sigma\varphi -\varphi(\hat{\sigma})\big)\dvol_{\Delta_n} \right)\\
&&+ \left(\frac{f_p(\Sd^d(\Delta_n))}{(n+1)!^d}-q_{p,n}\right)\frac{1}{(n+1)!^l}\sum_{\sigma\in \Sd^l(\Delta_n)^{[n]}} \varphi(\hat{\sigma})-\int_{\Delta_n}\varphi d\theta^l_p(d).\end{array}$

Now, since $\varphi$ is continuous, $\sup_{\sigma\in \Sd^l(\Delta_n)^{[n]}}(\sup_{\sigma}|\varphi-\varphi(\hat{\sigma})|)$ converges to 0 as $l$ grows to $+\infty$, while $\frac{1}{(n+1)!^l}|\sum_{\sigma\in \Sd^l(\Delta_n)^{[n]}}\varphi(\hat{\sigma})|$ remains bounded by $\sup_{\Delta_n}|\varphi|.$ Likewise by Theorem~\ref{Thm_DPS}, $\frac{f_p^d(\Delta_n)}{(n+1)!^d}$ converges to $q_{p,n}$ as $d$ grows to $+\infty$, while by Lemma~\ref{Lem_theta}, 
$\int_{\Delta_n} 1d\theta_p^l(d)$ converges to 0.
By letting $d$ grow to $+\infty$ and then $l$ grow to $+\infty$, we deduce that  $R_{l,d}$ can be as small as we want for $l,d$ large enough. This proves the result for $K=\Delta_n$.

Now, if $K$ is a locally finite  $n$-dimensional simplicial complex, we deduce the result by summing over all $n$-dimensional simplices of $K$, since from Theorem~\ref{Thm_DPS}, the measure of the $(n-1)$-skeleton  of $K$ with respect to $\gamma_{p,K}^d$ converges to 0 as $d$ grows to $+\infty$.
\hfill$\square$\\

Note that by integration of the constant function 1, Theorem~\ref{Thm_Gammapd} implies that for a finite simplicial complex $K$, $\frac{f_p^d(K)}{(n+1)!^d}  \underset{d \to + \infty}\longrightarrow q_{p,n},$ recovering the first part of Theorem~\ref{Thm_DPS}.
Also, since $q_{n,n}=1$, it implies that $\gamma_{n,K}^d \underset{d \to + \infty}\longrightarrow \dvol_K.$
This actually quickly follows from Riemann integration, since for every $\varphi\in C^0(\Delta_n),$

$$\begin{array}{lcl} 
\int_{\Delta_n}\varphi\dvol_{\Delta_n}&=&\lim_{d\to +\infty}\frac{1}{(n+1)!^d}\sum\limits_{\sigma\in\Sd^d(\Delta_n)^{[n]}}\varphi(\hat{\sigma})\\
&&\\
&=&\lim_{d\to +\infty}\int_{\Delta_n}\varphi d\gamma_{n,\Delta_n}^d.
\end{array}$$

Let us give another point of view of this fact.
For every $\sigma\in \Sd(\Delta_n)^{[n]},$ let us choose once for all a simplicial isomorphism $f_\sigma :\Delta_n\to \sigma.$
Let us  then consider  the product space $\Omega=Map(\N^*, \Sd(\Delta_n)^{[n]})=(\Sd(\Delta_n)^{[n]})^{\N^*}$ of countably many copies of $\Sd(\Delta_n)^{[n]}$ and equip it with the product measure $\omega$, where each copy of $\Sd(\Delta_n)^{[n]}$ is equipped with the counting measure  $\frac{1}{(n+1)!}\sum_{\sigma\in \Sd(\Delta_n)^{[n]}}\delta_\sigma.$ It is a Radon measure with respect to the product topology on $\Omega.$ We then set 

$$\begin{array}{rlcl}
\Phi :& \Omega\times \Delta_n& \to& \Delta_n\\
&((\sigma_i)_{i\in{\N^*}}, x)&\mapsto& \lim_{d\to +\infty} f_{\sigma_1}\circ\ldots\circ f_{\sigma_d}(x).
\end{array}$$

\begin{thm}\label{Thm_Phi}
The map $\Phi$ is well defined, continuous, surjective  and contracts the second factor $\Delta_n.$ Moreover,
$\dvol_{\Delta_n}=\Phi_*(\omega\times dvol_{\Delta_n})=\Phi_*(\omega\times \delta_{\hat{\Delta}_n})=\lim_{d\to +\infty}\gamma_{n, \Delta_n}^d.$ 
\end{thm}
(This result may be compared to the general Borel isomorphism theorem.) 

For every $d\geq 1$, let us set 

$$\begin{array}{rlcl}
\Phi_d :& \Omega\times \Delta_n& \to& \Delta_n\\
&((\sigma_i)_{i\in{\N^*}}, x)&\mapsto&  f_{\sigma_1}\circ\ldots\circ f_{\sigma_d}(x).
\end{array}$$

\begin{pf}
For every $(\sigma_i)_{i\in \N^*}\in \Omega$, the sequence of compact subsets $\Image(f_{\sigma_1}\circ\ldots\circ f_{\sigma_d})$  decreases as $d$ grows to $+\infty.$ These subsets are $n$-simplices of the barycentric subdivision $\Sd^d(\Delta_n)$ so that  their diameters converge to zero. We deduce the first part of Theorem~\ref{Thm_Phi}. 
Since $\Phi$  contracts the  second factor and is measurable, the push forward $\Phi_*(\omega\times \mu)$ does not depend on the probability measure $\mu$ on $\Delta_n$. In particular, $\Phi_*(\omega\times \dvol_{\Delta_n})=\Phi_*(\omega\times \delta_{\hat{\Delta}_n}).$
Now, we have by definition $(\Phi_d)_*(\omega\times \dvol_{\Delta_n})=\frac{1}{(n+1)!^d}\sum_{\tau\in \Sd^d(\Delta_n)}(f_\tau)_*(\dvol_{\Delta_n})$, where $f_\tau$ is the corresponding  simplicial isomorphism $f_{\sigma_1}\circ\ldots\circ f_{\sigma_d}$ between $\Delta_n$ and $\tau,$ so that $(\Phi_d)*(\omega\times \dvol_{\Delta_n})= \dvol_{\Delta_n}$ for every $d$ since $(f_\tau)_*(\dvol_{\Delta_n})=(n+1)!^d\dvol_{\Delta_n}|_\tau.$
Likewise, $(\Phi_d)_*(\omega\times \delta_{\hat{\Delta}_n})=\frac{1}{(n+1)!^d}\sum_{\tau\in \Sd^d(\Delta_n)}(f_\tau)_*(\delta_{\hat{\Delta}_n})=\gamma_{n,\Delta_n}^d$ by definition. Since the sequence $(\Phi_d)_{d\in \N^*}$ of continuous maps converge pointwise to $\Phi$, we deduce from Lebesgue's dominated convergence theorem that  for every probability measure $\mu$ on $\Delta_n,$ the sequence $(\Phi_d)_*(\omega\times \mu)$ weakly converges to $\Phi_*(\omega\times \mu).$
 \end{pf}

Recall that by definition, the Dirac measure $\delta_{\hat{\Delta}_n}$ in Theorem~\ref{Thm_Phi} coincides with  the measure $\gamma_{n,\Delta_n}.$ For $p<n$, we get 

\begin{thm} \label{Thm_fpDeltan}
For every $p\in \{0,\ldots,n\},$
$$f_p(\Delta_n)\dvol_{\Delta_n}=\Phi_*(\omega\times \gamma_{p,\Delta_n})=\lim_{d\to +\infty}f_{n-p-1}\big(\Lk(\sigma,\Sd^d(\Delta_n))\big)d\gamma^d_{p,\Delta_n}(\sigma).$$
\end{thm}

Recall that $f_p(\Delta_n)=\binom{n+1}{p+1}$ and that  by definition $f_{-1}(\Lk(\sigma, \Sd^d(\Delta_n)))=1$.

\begin{pf}
From Theorem~\ref{Thm_Phi}, $\Phi$ contracts the second factor. Since the mass of $\gamma_{p, \Delta_n}$ equals $f_p(\Delta_n)$ by definition, we deduce the first equality. Now, as in the proof of Theorem~\ref{Thm_Phi}, we deduce from Lebesgue's dominated convergence theorem that the sequence $(\Phi_d)_*(\omega\times \gamma_{p,\Delta_n})$ weakly converges to $\Phi_*(\omega\times \gamma_{p,\Delta_n}).$  It remains thus to compute $(\Phi_d)_*(\omega\times \gamma_{p,\Delta_n}).$ By definition   $(\Phi_d)_*(\omega\times \gamma_{p,\Delta_n})=\frac{1}{(n+1)!^d}\sum_{\tau\in \Sd^d(\Delta_n)^{[n]}}(f_\tau)_*(\gamma_{p,\Delta_n}),$ where $f_\tau$ is the corresponding simplicial isomorphism $f_{\sigma_1}\circ\ldots\circ f_{\sigma_d}$ between $\Delta_n$ and $\tau$. In this sum, we see that each $p$-simplex of $\Sd^d(\Delta_n)$ receives as many Dirac measures as the number of $n$-simplices adjacent to it. The number of  $n$-simplices adjacent to $\sigma\in \Sd^d(\Delta_n)^{[p]}$   is by definition $f_{n-p-1}(\Lk(\sigma, \Sd^d(\Delta_n))).$ We deduce 

$$\begin{array}{rcl}(\Phi_d)_*(\omega\times \gamma_{p,\Delta_n})&=&\frac{1}{(n+1)!^d}\sum\limits_{\sigma\in \Sd^d(\Delta_n)^{[p]}}f_{n-p-1}\big(\Lk(\sigma, \Sd^d(\Delta_n))\big)\delta_{\hat{\sigma}}\\
&&\\
&=&f_{n-p-1}\big(\Lk(\sigma, \Sd^d(\Delta_n))\big)d\gamma^d_{p,\Delta_n}(\sigma).\end{array}$$ 
\end{pf}

\begin{cor}\label{Cor_fpDeltan}
For every  $n$-dimensional locally finite simplicial complex $K$ and every $p\in \{0,\ldots, n\}$, the measure $f_{n-p-1}\big(\Lk(\sigma, \Sd^d(K))\big)d\gamma^d_{p,K}(\sigma)$ weakly converges to $f_p(\Delta_n)\dvol_K$ as $d$ grows to $+\infty.$
\end{cor}

\begin{pf} By definition $$\gamma^d_{p,K}=\sum_{\sigma\in K^{[n]}}\gamma^d_{p,\sigma} -\sum_{\tau \in K^{(n-1)}}\left( f_{n-\dim\tau-1}(\Lk(\tau, K))-1\right)\left(\frac{(\dim \tau+1)!}{(n+1)!}\right)^d\gamma_{p,\tau}^d$$ since  for every $\tau\in K^{(n-1)}$  and every $\sigma\in K^{[n]}$ such that $\tau<\sigma$, $\gamma_{p,\sigma}^d|_{\tau}= \left(\frac{(\dim \tau+1)!}{(n+1)!}\right)^d\gamma_{p,\tau}^d$ by definition and $\tau$ is a face of exactly $f_{n-\dim\tau-1}(\Lk(\tau, K))$ such $\sigma'$s. The result thus follows from Theorem~\ref{Thm_Gammapd} and Theorem~\ref{Thm_fpDeltan}. \end{pf}

\section{Limit density of links in a simplicial complex}\label{Sect_Links}
Corollary~\ref{Cor_fpDeltan} computes the limit density as $d$ grows to $+\infty$ of the top face numbers of the links of $p$-dimensional simplices in $\Sd^d(K), p\in \{0,\ldots, n\}$. We are going now to extend this result to all the face numbers of these links.

\begin{thm}\label{Thm_qLk}
For every $n$-dimensional  locally finite simplicial complex $K$ and every $0\leq p<n,$ the measure 
$q_{\Lk(\sigma,\Sd^d(K))}(T)d\gamma_{p,K}^d(\sigma)$  (with value in $\R_{n-p-1}[T]$) weakly converges to $\left(\sum_{l=0}^{n-p-1}q_{p+l+1,n} f_p(\Delta_{p+l+1})T^l\right)\dvol_K$ as $d$ grows to $+\infty.$
\end{thm}
\begin{pf} Let $\varphi\in C_c^0(|K|)$ be a continuous function with compact support on $|K|$. For every $0\leq l\leq n-p-1$, let us introduce the set 
\begin{align}
\mathcal{I}_l=\{(\sigma, \tau)\in \Sd^d(K)^{[p]}\times \Sd^d(K)^{[p+l+1]}| \sigma< \tau\}.\label{Defn_Il}\end{align}
It is equipped with the projection $p_1: (\sigma,\tau)\in \mathcal{I}_l\mapsto \sigma\in \Sd^d(K)^{[p]}$ and  $p_2: (\sigma,\tau)\in \mathcal{I}_l\mapsto \tau\in \Sd^d(K)^{[p+l+1]}$. We observe that for every  $\sigma\in \Sd^d(K)^{[p]}$,
$\# p_1^{-1}(\sigma)=f_l(\Lk(\sigma, \Sd^d(K)))$ while for every $\tau\in  \Sd^d(K)^{[p+l+1]}$, $p_2^{-1}(\tau)$ is in bijection with $\tau^{[p]}$ (given by $p_1$). Let us set 

\begin{align}
\varphi_1: (\sigma, \tau)\in \mathcal{I}_l\mapsto \varphi(\hat{\sigma})\in \R;\,\,\,
\varphi_2: (\sigma, \tau)\in \mathcal{I}_l\mapsto \varphi(\hat{\tau})\in \R;\,\,\, \gamma_l=\frac{1}{(n+1)!^d}\sum_{(\sigma,\tau)\in \mathcal{I}_l}\delta_{(\sigma,\tau)}.\label{Defn_gammal}\end{align}

Then, we deduce 
$$\begin{array}{rcl}
\int_{K}\varphi f_l\big(\Lk(\sigma, \Sd^d(K))\big)d\gamma_{p,K}^d(\sigma)&=&\int_{\mathcal{I}_l}\varphi_1d\gamma_l\\
&=&\int_{\mathcal{I}_l}\varphi_2d\gamma_l+\int_{\mathcal{I}_l}(\varphi_1-\varphi_2)d\gamma_l\\
&=&\int_{\Sd^d(K)^{[p+l+1]}}(p_2)_*(\varphi_2d\gamma_l)+ \int_{\mathcal{I}_l}(\varphi_1-\varphi_2)d\gamma_l\\
&=&\int_{K} \varphi f_p(\tau)d\gamma^d_{p+l+1,K}(\tau)+ \int_{\mathcal{I}_l}(\varphi_1-\varphi_2)d\gamma_l
\end{array}
$$ 
From Theorem~\ref{Thm_Gammapd}, the first term $\int_{K} \varphi f_p(\tau)d\gamma^d_{p+l+1,K}(\tau)$ in the right hand side converges to $q_{p+l+1,n}f_p(\Delta_{p+l+1})\int_K\varphi\dvol_K$ as $d$ grows to $+\infty$ while the second term $\int_{\mathcal{I}_l}(\varphi_1-\varphi_2)d\gamma_l
$ converges to zero.
Indeed, $\varphi$ is continuous with compact support and the diameter of $\tau\in \Sd^d(K)^{[p+l+1]}$ uniformly converges to zero on this compact  subset as $d$ grows to $+\infty.$  Thus, the suppremum of $(\varphi_1-\varphi_2)$ converges to zero as $d$ grows $+\infty.$
On the other hand, the total mass of $\gamma_l$ remains bounded, since $$ \int_{\mathcal{I}_l}1\gamma_l=\int_{\Sd^d(K)^{[p+l+1]}}(p_2)_*(d\gamma_l)=f_p(\Delta_{p+l+1})\int_K\gamma_{p+l+1,K}^d$$
and the latter is bounded from Theorem~\ref{Thm_Gammapd}. The result follows by definition of $q_{\Lk(\sigma, \Sd^{d}(K))}(T)$.
\end{pf}
Note that the $(n-1)$-skeleton of $K$ has vanishing measure with respect to $\dvol_K$ while for every $\sigma\in \Sd^d(K)^{[p]}$ interior to an $n$-simplex, its link is a homology $(n-p-1)$-sphere (Theorem~63.2 of~\cite{M}). After evaluation at $T=-1 $ and integration  of the constant function 1,    Theorem~\ref{Thm_qLk} thus provides the following asymptotic Dehn-Sommerville relations:

$$\sum_{l=p}^nq_{l,n}\binom{l+1}{p+1}(-1)^{n+l}=q_{p,n}.$$

Now,  recall that the dual  block $D(\sigma)$ of a simplex $\sigma\in K$ is the union of all open simplices $[\hat{\sigma}_0\ldots, \hat{\sigma}_p]$ of $\Sd(K)$ such that $\sigma_0=\sigma$, see~\cite{M}.
The closure $\overline{D}(\sigma)$ of $D(\sigma)$ is called \emph{closed block dual} to $\sigma$ and  following~\cite{M}  we set $\dot{D}(\sigma)=\overline{D}(\sigma)\setminus D(\sigma).$
Then, we get the following.

\begin{thm}\label{Thm_qDsigma}
For every $n$-dimensional locally finite simplicial complex $K$ and every $0\leq p\leq n,$
 the measure $q_{D(\sigma)}(T)d\gamma^d_{p,K}(\sigma)$ weakly converges to  $\sum_{l=0}^{n-p}\left(\sum_{h=l}^{n-p}q_{p+h,n}f_p{(\Delta_{p+h})}\lambda_{h,l}\right)T^l \dvol_K$ as $d$ grows to $+\infty.$
\end{thm}

\begin{pf} By definition,  the dual block $D(\sigma)$ has only one face in dimension 0, namely $\hat{\sigma}$, so that for the coefficient $l=0$, the result follows from Theorem~\ref{Thm_Gammapd}. Let us now assume  that $0<l\leq n-p$ and choose
$\varphi\in C_c^0(|K|).$ We set 

$$\mathcal{J}_l=\{(\sigma, \theta)\in \Sd^d(K)^{[p]}\times \Sd^{d+1}(K)^{[l-1]} | \theta\in \dot{D}(\sigma)\}.$$
Let $p_1 : (\sigma, \theta) \in \mathcal{J}_l\mapsto \sigma\in \Sd^d(K)^{[p]}$. Then, for every $\sigma\in \Sd^d(K)^{[p]}$,
$\#p_1^{-1}(\sigma)=f_l(D(\sigma))$, since $p_1^{-1}(\sigma)$
is in bijection with $\dot{D}(\sigma)$ and by taking the cone over $\hat{\sigma}$ we get an isomorphism $\tau\in \dot{D}(\sigma)\mapsto \hat{\sigma}\ast\tau \in D(\sigma)\setminus\hat{\sigma}$ where $\ast$ denotes the join operation. (Recall that if $\tau=[e_0,\ldots, e_k]$ the join $\hat{\sigma}\ast\tau$ is $[\hat{\sigma}, e_0,\ldots, e_k].$)

 Likewise by definition, every simplex $\theta\in \dot{D}(\sigma)^{[l-1]}$ reads $\theta=[\hat{\tau}_0,\ldots, \hat{\tau}_{l-1}]$ where $\sigma<\tau_0<\ldots<\tau_{l-1}$ are simplices of $\Sd^d(K)$ (see Theorem~64.1 of~\cite{M}). We deduce a map 
 
$$\begin{array}{rccl}
\pi: &\mathcal{J}_l&\to &\bigsqcup\limits_{h=l-1}^{n-p-1}\mathcal{I}_h\\
&(\sigma, [\hat{\tau}_0,\ldots, \hat{\tau}_{l-1}])&\mapsto& (\sigma, \tau_{l-1})
\end{array}
$$
where $\mathcal{I}_h$ is the set defined  in (\ref{Defn_Il}).

We then set $p_2 :(\sigma,\tau)\in \bigsqcup_{h=l-1}^{n-p-1}\mathcal{I}_h \mapsto \tau\in \Sd^d(K)\setminus \Sd^d(K)^{(p+l-1)}.$ As in the proof of Theorem~\ref{Thm_qLk}, for every $\tau\in \Sd^d(K)\setminus \Sd^d(K)^{(p+l-1)},$ 
$p_2^{-1}(\tau)$ is in bijection with $\tau^{[p]}$ and $\pi^{-1}(\sigma,\tau)$ with the set of interior $(l-1)$-dimensional simplices of $\Sd(\Lk(\sigma, \tau))$,  so that  $\#\pi^{-1}((\sigma, \tau))=\lambda_{h+1,l}$ if $\dim\tau=p+h+1.$
Let us set $\tilde{\varphi}_1 : (\sigma, \tau)\in\mathcal{J}_l\mapsto \varphi(\hat{\sigma})\in \R$ and $\tilde{\gamma}_l=\frac{1}{(n+1)!^d} \sum_{(\sigma, \theta)\in \mathcal{J}_l}\delta_{(\sigma, \theta)}$. Then, we deduce 

$$\begin{array}{rcl}
\int_K \varphi f_l(D(\sigma))d\gamma_{p,K}^d(\sigma)&=&\int_{\mathcal{J}_l} \tilde{\varphi}_1 d\tilde{\gamma}_l\\
 &=&\sum_{h=l-1}^{n-p-1}\lambda_{h+1, l}\int_{\mathcal{I}_h}\varphi_1d\gamma_h 
 \end{array}
 $$
by pushing forward $\tilde{\varphi}_1d\tilde{\gamma}_l$ onto $\bigsqcup_{h=l-1}^{n-p-1}\mathcal{I}_h$ with $\pi$, where $\varphi_1$ and $\gamma_h$
are defined by (\ref{Defn_gammal}).

Now, we have established in the proof of Theorem~\ref{Thm_qLk} that as $d$ grows to $+\infty$, $\int_{\mathcal{I}_h}\varphi_1d\gamma_h$ converges to  $f_p(\Delta_{p+h+1})q_{p+h+1}\int_{K}\varphi \dvol_K$. We deduce that  $f_l(D(\sigma))d\gamma_{p,K}^d(\sigma)$ weakly converges to $\left(\sum_{h=l}^{n-p}\lambda_{h,l}f_p(\Delta_{p+h}) q_{p+h,n}\right)\dvol_K$. Hence the result.
\end{pf}

\begin{rem} In~\cite{SW2}, we study the expected topology of a random subcomplex in a finite simplicial complex $K$ and its barycentric subdivisions.
The Betti numbers  of such a subcomplex turn out to be asymptotically controlled by the measure given by Theorem~\ref{Thm_qDsigma}.
\end{rem}

 Let us now finally observe that Theorem~\ref{Thm_qDsigma} provides a geometric proof of the following (compare Theorem~A of~\cite{DPS}).

\begin{cor}\label{Cor_Eigenvector}
The vector $(q_{p,n})_{0\leq p\leq n}$ is the eigenvector of $\Lambda_n^t$ associated to the eigenvalue $(n+1)!$, normalized by the relation $q_{n,n}=1$.
\end{cor}

\begin{pf}
By Theorem~64.1 of~\cite{M}, we know that the dual blocks of a complex $K$ are disjoint and that their union is $|K|$. We deduce that for every $d\in\N^*$, 

$$\frac{1}{(n+1)!^d}q_{\Sd^{d+1}(\Delta_n)}(T)=\sum_{p=0}^{n}\int_{\Delta_n}q_{D(\sigma)}(T)d\gamma^d_{p,\Delta_n}(\sigma).$$
By letting $d$ grow to $+\infty$, we now deduce from Theorem~\ref{Thm_qDsigma}, applied to $K=\Delta_n$ and  after integration of 1, that

$$\begin{array}{rcl}

(n+1)!\sum\limits_{p=0}^{n}q_{p,n}T^p&=&\sum\limits_{p=0}^{n}\left( \sum\limits_{l=0}^{n-p}T^l \sum\limits_{h=p+l}^{n}q_{h,n}f_p(\Delta_h)\lambda_{h-p,l}\right)\\
&=&\sum\limits_{l=0}^{n}T^l \left( \sum\limits_{h=l}^{n}q_{h,n} \sum\limits_{p=0}^{h-l}f_p(\Delta_h)\lambda_{h-p,l}\right)\end{array}$$

Now, $\sum\limits_{p=0}^{h-l}f_p(\Delta_h)\lambda_{h-p,l}=\sum\limits_{p=l}^{h}\binom{h+1}{p}\lambda_{p,l}=\lambda_{h+1,l+1}$ from Lemma~\ref{Lem_lambdaij}. Hence, for every $p\in\{0,\ldots,n\}$, $(n+1)!q_{p,n}=\sum\limits_{h=l}^{n}q_{h,n}\lambda_{h+1,l+1}.$
\end{pf}
\bibliography{Random1}

\begin{thebibliography}{1}

\bibitem{A}
T.~Akita.
\newblock A formula for the {E}uler characteristics of even dimensional
  triangulated manifolds.
\newblock {\em Proc. Amer. Math. Soc.}, 136(7):2571--2573, 2008.

\bibitem{BW}
F.~Brenti and V.~Welker.
\newblock {$f$}-vectors of barycentric subdivisions.
\newblock {\em Math. Z.}, 259(4):849--865, 2008.

\bibitem{DPS}
E.~Delucchi, A.~Pixton, and L.~Sabalka.
\newblock Face vectors of subdivided simplicial complexes.
\newblock {\em Discrete Math.}, 312(2):248--257, 2012.

\bibitem{Klain}
S.~Klain.
\newblock Dehn-sommerville relations for triangulated manifolds.
\newblock {\em unpublished manuscript available at
  http://faculty.uml.edu/dklain/ds.pdf.}

\bibitem{Klee}
V.~Klee.
\newblock A combinatorial analogue of {P}oincar\'e's duality theorem.
\newblock {\em Canad. J. Math.}, 16:517--531, 1964.

\bibitem{Mac}
I.~G. Macdonald.
\newblock Polynomials associated with finite cell-complexes.
\newblock {\em J. London Math. Soc. (2)}, 4:181--192, 1971.

\bibitem{M}
J.~R. Munkres.
\newblock {\em Elements of algebraic topology}.
\newblock Addison-Wesley Publishing Company, Menlo Park, CA, 1984.

\bibitem{SW2}
N.~Salepci and J.-Y. Welschinger.
\newblock Asymptotic topology of random subcomplexes in a finite simplicial
  complex.
\newblock {\em In preparation}, 2017.

\end{thebibliography}
\bibliographystyle{abbrv}

Univ Lyon, Universit\'e Claude Bernard Lyon 1, CNRS UMR 5208, Institut Camille Jordan, 43 blvd. du 11 novembre 1918, F-69622 Villeurbanne cedex, France

{salepci@math.univ-lyon1.fr, welschinger@math.univ-lyon1.fr.}
\end{document}